\documentclass[a4paper]{amsart}  

\usepackage{amssymb} \usepackage{amscd} \usepackage{times}

\def\zbb{\mathbb{Z}}  
  
  \def\phi{\varphi}
 \def\p1{{\mathbb{P}^1_\zbb}}

\newtheorem{Theorem}{\quad Theorem}[section]

\newtheorem{Lemma}[Theorem]{\quad Lemma}

\newtheorem{Remark}[Theorem]{\quad Remark}

\begin{document}

\title{sup $ \times $ inf  inequalities for the scalar curvature equation in dimensions 4 and 5.}

\author{Samy Skander Bahoura}

\address{Departement de Mathematiques, Universite Pierre et Marie Curie, 2 place Jussieu, 75005, Paris, France.}
              
\email{samybahoura@yahoo.fr} 

\date{}

\maketitle

\begin{abstract}

We consider the following problem on bounded open set  $ \Omega $ of $ {\mathbb R}^n $:

$$ \left \{ \begin {split} 
     -\Delta u & = Vu^{\frac{n+2}{n-2}} \,\, &&\text{in} \!\!&&\Omega \subset {\mathbb R}^n, \,\, n=4, 5, \\
                  u & > 0  \,\,                                 && \text{in} \!\!&&\Omega.               
\end {split}\right. $$


We assume that : 

$$ V \in C^{1, \beta}(\Omega), \,\,  0 < \beta \leq 1 $$ 

$$  0 < a \leq V \leq b  < + \infty, $$

$$ |\nabla V| \leq A \,\,\text{in} \,\,\Omega, $$

$$ |\nabla^{1+ \beta} V| \leq B \,\, \text{in} \,\, \Omega.$$

then, we have a $  \sup  \times  \inf $ inequality for the solutions of the previous equation, namely:

$$ (\sup_K u)^{\beta} \times  \inf_{\Omega} u \leq c=c(a, b, A, B, \beta, K, \Omega), \,\, \text{for} \,\, n = 4, $$

and,

$$ (\sup_K u)^{1/3} \times  \inf_{\Omega} u \leq c=c(a, b, A, B, K, \Omega), \,\, \text{for} \,\, n = 5, \,\, \text{and} \,\,\beta=1. $$

\end{abstract}

\section{Introduction and Main Result}

We work on $ \Omega \subset  \subset {\mathbb R}^4 $ and we consider the following equation:

\begin{displaymath}\left \{ \begin {split} 
     -\Delta u & = Vu^{\frac{n+2}{n-2}} \,\, &&\text{in} \!\!&&\Omega \subset {\mathbb R}^n, \,\, n=4, 5, \\
                  u & > 0  \,\,                                 && \text{in} \!\!&&\Omega.               
\end {split}\right.\qquad (E) \end{displaymath}

with,

$$   \begin {cases}

 V \in C^{1, \beta}(\Omega),  \\
 
 0 < a \leq V \leq b  < + \infty \,\,\text{in} \,\,\Omega,  \\
 
 |\nabla V| \leq A \,\,\text{in} \,\,\Omega, \\
  
 |\nabla^{1+ \beta}V| \leq B \,\,\text{in} \,\,\Omega.

\end {cases} \qquad (C_{\beta }) $$

Without loss of genarality, we suppose $ \Omega = B_1(0) $ the unit ball of $ {\mathbb R}^n $.

\smallskip

The corresponding equation in two dimensions on open set $ \Omega $ of $ {\mathbb R}^2 $, is:

$$ -\Delta u=V(x)e^u, \qquad (E') $$

The equation $ (E') $ was studied by many authors and we can find very important result about a priori estimates in [8], [9], [12], [16], and [19]. In particular in [9] we have the following interior estimate:

$$  \sup_K u  \leq c=c(\inf_{\Omega} V, ||V||_{L^{\infty}(\Omega)}, \inf_{\Omega} u, K, \Omega). $$

And, precisely, in [8], [12], [16], and [19], we have:

$$ C \sup_K u + \inf_{\Omega} u \leq c=c(\inf_{\Omega} V, ||V||_{L^{\infty}(\Omega)}, K, \Omega), $$

and,

$$ \sup_K u + \inf_{\Omega} u \leq c=c(\inf_{\Omega} V, ||V||_{C^{\alpha}(\Omega)}, K, \Omega). $$

where $ K $ is a compact subset of $ \Omega $, $ C  $ is a positive constant which depends on $\dfrac{\inf_{\Omega} V}{\sup_{\Omega} V} $, and,  $ \alpha \in (0, 1] $.

For $ n  \geq 3 $ we have the following general equation on a riemannian manifold:

$$ -\Delta u+hu=V(x) u^{\frac{n+2}{n-2}},\,\, u >0 . \qquad (E_n) $$

Where $ h, V $ are two continuous functions. In the case  $ c_n h= R_g $ the scalar curvature, we call $ V $ the prescribed scalar curvature. Here $ c_n $ is a universal constant.

\bigskip

The equation $ (E_n) $ was studied a lot, when $ M =\Omega \subset {\mathbb R}^n $ or $ M={\mathbb S}_n $ see for example, [2-4], [11], [15]. In this case we have a $ \sup \times \inf $ inequality.

\smallskip

In the case $ V\equiv 1 $ and $ M $ compact, the equation $ (E_n) $ is Yamabe equation. T.Aubin and R.Schoen proved the existence of solution in this case, see for example [1] and [14] for a complete and detailed summary.

\bigskip

When $ M $ is a compact Riemannian manifold, there exist some compactness result for equation  $ (E_n) $ see [18]. Li and Zhu see [18], proved that the energy is bounded and if we suppose $ M $ not diffeormorfic to the three sphere, the solutions are uniformly bounded. To have this result they use the positive mass theorem.

\bigskip

Now, if we suppose $ M $ Riemannian manifold (not necessarily compact) and $ V\equiv 1 $, Li and Zhang [17] proved that the product $ \sup \times \inf $ is bounded. On other handm see [3], [5] and [6] for other Harnack type inequalities, and, see [3] and [7] about some caracterisation of the solutions of this equation $ (E_n) $ in this case ($ V\equiv 1 $).

Here we extend a result of [11] on an open set of  $ {\mathbb R}^n, n= 4, 5 $. In fact we consider the prescribed scalar curvature equation on an open set of  $ {\mathbb R}^n, n= 4, 5 $, and, we prove a $  \sup  \times \inf $ inequality on compact set of the domain when the derivative of the prescribed scalar curvature is $\beta $-holderian, $ \beta >0 $.

Our proof is an extension of  Chen-Lin result in dimension 4 and 5, see [11] , and,  the moving-plane method is used to have this estimate. We refer to Gidas-Ni-Nirenberg for the  moving-plane method, see  [13]. Also, we can see in [10], one of the application of this method.
 
 \bigskip 
 
We have the following result in dimension 4, which is the consequence of the work of  Chen-Lin.

\medskip

{\bf Theorem A}. For all $ a,b, m, A, B >0 $, and for all compact $ K $ of $ \Omega $, there exists a positive constant $ c=c(a, b, A, B, K, \Omega) $ such that:

$$ \sup_K u \times  \inf _{\Omega }  \leq c, $$

where $ u $ is solution of $ (E) $ with $ V $, $ C^{2}$ satisfying $ (C_{\beta}) $ for $ \beta = 1 $.
 
\bigskip

Here, we give an inequality of type $ \sup \times \inf $ for the equation $ (E) $  in dimension 4 and with general conditions  on  the prescribed scalar curvature, exactly we take a $ C^{1, \beta} $ condition. In fact we extend the result of Chen-Lin in dimension 4.

\medskip

Here we prove:

\medskip

\begin{Theorem}. For all $ a,b, A, B >0 $, $ 1 \geq \beta >0 $, and for all compact $ K $ of $ \Omega $, there exists a positive constant $ c=c(a, b, A, B, \beta,  K, \Omega) $ such that:

$$ (\sup_K u)^{\beta} \times  \inf_{\Omega} u \leq c, $$

where $ u $ is solution of $ (E) $ with $ V$ satisfying $ (C_{\beta}) $. 

\end{Theorem}

\bigskip

We have the following result in dimension 5, which is the consequence of the work of  Chen-Lin.

\medskip

{\bf Theorem B}. For all $ a,b, m, A, B >0 $, and for all compact $ K $ of $ \Omega $, there exists a positive constant $ c=c(a, b, m, A, B, K, \Omega) $ such that:

$$ \sup_K u \leq c,  \,\,  {\rm if }  \,\, \inf_{\Omega} u \geq m, $$

where $ u $ is solution of $ (E) $ with $ V$ satisfying $ (C_{\beta})= (C_1) $ for $ \beta=1 $.
 
\bigskip

Here, we give an inequality of type $ \sup \times \inf $ for the equation $ (E) $  in dimension 5 and with general conditions  on  the prescribed scalar curvature, exactly we take a $ C^2 $ condition ($ \beta = 1$ in $ (C_{\beta }) $). In fact we extend the result of Chen-Lin in dimension 5.

\medskip

Here we prove:

\medskip

\begin{Theorem}. For all $ a,b, A, B >0 $, and for all compact $ K $ of $ \Omega $, there exists a positive constant $ c=c(a, b, A, B, K, \Omega) $ such that:

$$ (\sup_K u)^{1/3} \times  \inf_{\Omega} u \leq c, $$

where $ u $ is solution of $ (E) $ with $ V$ satisfying $ (C_{\beta}) $ for $ \beta=1 $.

\end{Theorem}

\bigskip

\section{The method of moving-plane.}

In this section we will formulate a modified version of the method of moving-plane for use later. Let $  \Omega $ an open set and $  \Omega^c $ the complement of $  \Omega $. We consider a solution $ u $ of the following equation:

$$   \begin {cases}  
\Delta u + f(x, u)= 0,  \\
          \,\,\, u >0,
   \end{cases}  \qquad (E'') $$

where $ f(x, u) $ is nonegative, Holder continuous in $ x $, $ C^1 $ in $ u $, and defined on $  \bar\Omega  \times (0, +  \infty ) $. Let $ e $ be a unit vector in $ {\mathbb R}^n $. For $  \lambda <0 $, we let 

$ T_{\lambda} = \{ x \in {\mathbb R}^n, \langle x,e\rangle = \lambda \} $, $ \Sigma_{\lambda} = \{ x \in {\mathbb R}^n, \langle x,e\rangle > \lambda \} $, and $ x^{\lambda} = x+ (2 \lambda-2 \langle x,e\rangle)e $ to denote the reflexion point of $ x $ with respect to $ T_{\lambda} $, where $ \langle.,.\rangle $ is the standard inner product of $ {\mathbb R}^n $. Define:

$$ \lambda_1 \equiv  \sup \{ \lambda < 0,  \Omega^c \subset \Sigma_{\lambda}  \}, $$

$ \Sigma'_{\lambda} = \Sigma_{\lambda} - \Omega^c  $ for $ \lambda \leq \lambda_1 $, and $ \bar \Sigma'_{\lambda}  $ the closure of  $ \Sigma'_{\lambda} $. Let $ u^{\lambda} (x)= u(x^{\lambda}) $ and $ w_{\lambda}(x)=u(x)-u^{\lambda}(x) $ for $ x \in \Sigma'_{\lambda} $. Then we have, for any arbitrary function $ b_{\lambda}(x) $,

$$ \Delta w_{\lambda}(x)+ b_{\lambda}(x) w_{\lambda}(x) = Q(x, b_{\lambda}(x)), $$

where,

$$  Q(x, b_{\lambda}(x))= f(x^{\lambda}, u^{\lambda})-f(x,u)+ b_{\lambda}(x) w_{\lambda}(x). $$

The hypothesis $ (*) $ is said to be satisfied if there are two families of functions $ b_{\lambda}(x) $ and $ h^{\lambda}(x)  $ defined in $  \Sigma'_{\lambda} $, for $ \lambda \in (-\infty, \lambda_1) $ such that, the following assertions holds:

  $$ 0 \leq  b_{\lambda}(x) \leq c(x)|x|^{-2}, $$
  
  where $ c(x) $ is independant of  $ \lambda $ and tends to zero as $ |x| $ tends to $ + \infty $,
  
  $$ h^{\lambda}(x) \in C^1(\Sigma_{\lambda} \cap \Omega), $$
  
 and satisfies:
 
  $$  \begin {cases}

     \Delta h^{\lambda}(x) \geq Q(x, b_{\lambda}(x)) \,\, \text{in} \,\,\Sigma_{\lambda} \cap \Omega \\
                  h^{\lambda}(x)  > 0 \,\, \text{in} \,\,  \Sigma_{\lambda} \cap \Omega
              
\end {cases} $$

in the distributional sense and,

$$ h^{\lambda}(x)= 0  \,\, \text{on} \,\,  T_{\lambda} \,\,\text{and} \,\, h^{\lambda}(x)= O(|x|^{-t_1}), $$

as $ |x| \to + \infty $ for some constant $ t_1 >0 $,

$$ h^{\lambda}(x) + \epsilon < w_{\lambda}(x), $$

in a neighborhood of $ \partial \Omega $, where $ \epsilon  $ is a positive constant independant of $ x $.

$$  \begin {cases}

     h^{\lambda}(x) \,\, \text{and} \,\, \nabla_x h^{\lambda} \,\, \text{are continuous with respect to both variables} \\
     
      \text{ $ x $ and $ \lambda $, and for any compact set of $ \Omega $,} \,\,   w_{\lambda}(x) > h^{\lambda}(x) \\
     
      \text{holds when $ -\lambda $ is sufficiently large.}             

\end {cases} $$

We have the following lemma:

\begin{Lemma} 

Let $ u $ be a solution of $ (E'') $. Suppose that $ u(x) \geq C >0 $ in a neighborhood of $ \partial \Omega $ and $ u(x)= O(|x|^{-t_2}) $ at $ +\infty $ for some positive $ t_2 $. Assume there exist $ b_{\lambda}(x) $ and $ h^{\lambda}(x) $ such that the hypothesis $ (*) $ is satisfied for $  \lambda  \leq  \lambda_1 $. Then $ w_{\lambda}(x) >0 $ in $ \Sigma'_{\lambda} $, and $ \langle \nabla u, e\rangle >0 $ on $ T_{\lambda} $ for $ \lambda \in (-\infty, \lambda_1) $.

\end{Lemma}

For the proof see Chen and Lin, [10].

\begin{Remark} 

If we know that $ w_{\lambda} - h^{\lambda} >0 $ for some $ \lambda= \lambda_0 < \lambda_1  $ and $ b_{\lambda} $ and $ h^{\lambda} $ satisfy the hypothesis $ (*) $ for $ \lambda_0 \leq \lambda \leq \lambda_1 $, then the conclusion of the lemma 2.1 holds.

\end{Remark}

\section{Proof of the result:} 

\underbar {Proof of the theorem 1, $ n=4 $ : }

\bigskip

To prove the theorem, we argue by contradiction and we assume that the $ (\sup )^{ \beta } \times \inf $ tends to infinity.

\bigskip

\underbar {Step 1: blow-up analysis }

\bigskip

We want to prove that:

$$  {\tilde R}^2 (\sup_{B_{\tilde R}(0)} u)^{ \beta } \times  \inf_{B_{3{\tilde R}}(0)} u \leq c= c(a, b, A, B, \beta ), $$

If it is not the case, we have:

$$ {\tilde R_i}^2 ( \sup_{B_{\tilde R_i}(0)} u_i)^{ \beta} \times  \inf_{B_{3 \tilde R_i}(0)} u_i = i^6 \to + \infty, $$

For positive solutions $ u_i >0 $ of the equation $ (E) $ and $ \tilde R_i \to 0. $

Thus,

$$ \dfrac{1}{i}{\tilde R_i}(\sup_{B_{\tilde R_i}(0)} u_i)^{(1+ \beta)/2} \to + \infty, $$

and,

$$  \dfrac{1}{i}{\tilde R_i} [\sup_{B_{\tilde R_i}(0)} u_i]^{(1+ \beta)/2 } \to + \infty, $$

Let $ a_i $ such that:

$$ u_i(a_i) = \max_{B_{\tilde R_i}(0)} u_i, $$

We set, 

$$ s_i(x)=(\tilde R_i-|x-a_i|)^{2/(1+ \beta )} u_i(x), $$

we have,

$$  s_i(\bar x_i)= \max_{B_{\tilde R_i}(a_i)} s_i \geq s_i(a_i)={\tilde R_i}^{2/ (1+ \beta )}\sup_{B_{R_i}(0)} u_i \to + \infty, $$

we set, 

$$ R_i= \dfrac{1}{2} (\tilde R_i-|\bar x_i-a_i|), $$

We have, for $ |x-\bar x_i|\leq \dfrac{R_i}{i} $, 

$$ \tilde R_i-|x-a_i|\geq \tilde R_i-|\bar x_i-a_i|-|x-a_i| \geq 2R_i-R_i=R_i$$

Thus,

$$ \dfrac{u_i(x)}{u_i(\bar x_i)} \leq \beta_i \leq 2^{2/(1+ \beta ) }. $$

with $ \beta_i \to 1. $

\bigskip

We set, 

 $$  M_i= u_i(\bar x_i),  \,\,\,v_i^*(y)=\dfrac{u_i(\bar x_i+M_i^{-1}y)}{u_i(\bar x_i)}, \,\,\, |y|\leq \dfrac{1}{i} R_iM_i^{(1+ \beta)/2 } = 2L_i. $$

And,

$$ \dfrac{1}{i^2}{R_i}^2 M_i^{ \beta } \times  \inf_{B_{3 \tilde R_i}(0)} u_i \to + \infty, $$

By the elliptic estimates, $ v_i^* $ converge on each compact set of $ {\mathbb R}^4 $ to a function $ U_0^* >0 $ solution of :

$$ \begin {cases}

   -\Delta U_0^* = V(0) {U_0^*}^{3} \,\, \text{in} \,\, {\mathbb R}^4, \\
                
                U_0^*(0)=1= \max_{{\mathbb R}^4} U_0^*.

\end {cases} $$

For simplicity, we assume that $ 0 < V(0)= n(n-2) =8 $. By a result of Caffarelli-Gidas-Spruck, see [10], we have:

$$ U_0^*(y)=(1+|y|^2)^{-1}. $$

We set,

$$ v_i(y)=v_i^*(y+e), $$

where $ v_i^* $ is the blow-up function. Then, $ v_i $ has a local maximum near $ -e $.

$$ U_0(y)=U_0^*(y+e). $$

We want to prove that:

$$ \min_{ \{0 \leq |y| \leq r  \} } v_i^*  \leq (1+ \epsilon) U_0^*(r). $$

for $ 0 \leq r \leq L_i $, with $  L_i= \dfrac{1}{2i}R_i M_i^{(1+\beta)/2 } $.

\bigskip

We assume that it is not true, then, there is a sequence of number $ r_i \in (0, L_i) $ and $ \epsilon >0 $, such that:

$$  \min_{ \{ 0 \leq |y| \leq r_i  \} } v_i^*  \geq (1+ \epsilon) U_0^*(r_i). $$

We have:

$$ r_i \to + \infty. $$

Thus , we have for $ r_i \in (0, L_i) $ :

$$  \min_{ \{ 0 \leq |y| \leq r_i  \} } v_i  \geq (1+ \epsilon) U_0(r_i). $$

Also, we can find a sequence of number $ l_i  \to + \infty $ such that:

$$ ||v_i^*-U_0||_{C^2(B_{l_i}(0))} \to 0. $$

Thus,

 $$ \min_{ \{ 0 \leq |y| \leq l_i  \} } v_i  \geq (1 - \epsilon/2) U_0(l_i). $$

\underbar {Step 2 : The Kelvin transform and the Moving-plane method }

\bigskip

\underbar {Step 2.1: a linear equation perturbed by a term, and, the auxiliary function } 

\bigskip

 {Step 2.1.1: $ D_i= |\nabla V_i(x_i)| \to 0.$} 

\bigskip

We have the same estimate as in the paper of Chen-Lin. We argue by contradiction. We consider $ r_i \in (0, L_i) $ where $ L_i $ is the number of the blow-up analysis.

$$  L_i= \dfrac{1}{2i}R_i M_i^{(1+ \beta )/2}. $$

 We use the assumption that the sup times inf is not bounded to prove $ w_{\lambda} >h_{\lambda} $ in $ \Sigma_{\lambda} =\{ y, y_1 > \lambda \}, $ and on the boundary.
 
 \bigskip

The function $ v_i $ has a local maximum near $ -e $ and converge to $ U_0(y)=U_0^*(y+e) $ on each compact set of $ {\mathbb R}^5 $. $U_0 $ has a maximum at $ -e $.

We argue by contradiction and we suppose that:

$$ D_i= |\nabla V_i(x_i)| \not  \to 0. $$

Then, without loss of generality we can assume that:

$$ \nabla V_i(x_i) \to e=(1,0,...0). $$

Where $ x_i $ is :

$$ x_i=\bar x_i+ M_i^{-1} e, $$

with $ \bar x_i  $ is the local maximum in the blow-up analysis.

\bigskip

As in the paper of Chen-Lin, we use the Kelvin transform twice and we set (we take the same notations):

$$ I_{\delta}(y)=\dfrac{\frac{|y|}{|y|^2}-\delta e}{\left (|\frac{|y|}{|y|^2}-\delta e| \right )^2}, $$

$$ v_i^{\delta}(y)=\dfrac {v_i(I_{\delta}(y))}{|y|^{n-2}|y-e/\delta|^{n-2}}, $$

and, 

$$  V_{\delta}(y) =V_i(x_i+ M_i^{-1}I_{\delta}(y)). $$

$$ U_{\delta}(y)=\dfrac {U_0(I_{\delta}(y))}{|y|^{n-2}|y-e/\delta|^{n-2}}. $$

Then, $ U_{\delta} $ has a local maximum near $ e_{\delta} \to -e $ when $ \delta \to 0 $. The function $ v_i^{\delta} $ has a local maximum near $ -e $.

We want to prove by the application of the maximum principle and the Hopf lemma that near $ e_{\delta} $ we have not a local maximum, which is a contradiction.

\bigskip

We set on $ \Sigma'_{\lambda} = \Sigma_{\lambda}  - \{ y, |y- \frac{e}{\delta}| \leq \frac{c_0}{r_i} \} \simeq \Sigma_{\lambda}  - \{ y, |I_{\delta}(y)| \geq  r_i \} $:

$$ h_{\lambda} (y)= - \int_{\Sigma_{\lambda}} G_{\lambda}(y, \eta) Q_{\lambda}(\eta) d\eta. $$
 
 with,

$$ Q_{\lambda}(\eta)=(V_{\delta}(\eta)-V_{\delta}(\eta^{\lambda})) (v_i^{\delta}(\eta^{\lambda}))^{3}. $$

And, by the same estimates, we have for $ \eta  \in A_1=\{ \eta, |\eta | \leq R=\epsilon_0/\delta \},$

$$ V_{\delta}(\eta)-V_{\delta}(\eta^{\lambda})  \geq M_i^{-1}(\eta_1-\lambda) + o(1)M_i^{-1}|\eta^{\lambda}| , $$

and, we have for $ \eta \in A_2 =\Sigma_{\lambda}-A_1 $:

$$ | V_{\delta}(\eta)-V_{\delta}(\eta^{\lambda}) |  \leq CM_i^{-1}(|I_{\delta}(\eta)| + |I_{\delta}(\eta^{\lambda})|), $$

And, we have for some $ \lambda_0 \leq -2 $ and $ C_0 >0 $:

$$ w_{\lambda }(y) = v_i^{\delta}(y)- v_i^{\delta}(y^{\lambda_0}) \geq C_0\dfrac{y_1-\lambda_0}{(1+|y|)^n}, $$

for $ y_1 > \lambda_0 $.

\bigskip

Because , by the maximum principle:

 $$ \min_{ \{ l_i \leq |I_{\delta}(y)| \leq  r_i  \} } v_i = \min  \{ \min_{ \{ |I_{\delta}(y)| = l_i  \} }, v_i  \min_{ \{  |I_{\delta}(y)| = r_i  \} } v_i  \}  \geq (1- \epsilon) U_{ \delta }(\dfrac{e}{\delta}) $$
 
 $$ \geq (1+c_1\delta - \epsilon) U_{ \delta }((\dfrac{e}{\delta})^{\lambda }) \geq (1+c_1\delta - 2\epsilon) v_i^{ \delta }(y^{\lambda }), $$
 
 and for $ |I_{\delta}(y)| \leq  l_i $ we use the $ C^2 $ convergence of $ v_i^{\delta} $  to $ U_{\delta} $. 
 
\bigskip

 Thus,
 
 $$ w_{\lambda }(y) >  2 \epsilon >0, $$ 
 
By the same estimates as in Chen-Lin paper (we apply the lemma 2.1 of the second section),  and by our hypothesis on $ v_i $, we have:

$$ 0 < h_{\lambda} (y) = O(1)M_i^{-1}(y_1- \lambda)(1+|y|)^{-n} < 2 \epsilon <  w_{\lambda}(y). $$

also, we have the same etimate on the boundary, $|I_{\delta}(\eta)| = r_i $ or  $ |y-e/\delta| =  c_2 r_i^{-1} $:

\bigskip

Step 2.1.1: $ |\nabla V_i(x_i)|^{1/ \beta} [u_i(x_i)] \leq C $

\bigskip

Here, also, we argue by contradiction. We use the same computation as in Chen-Lin paper,  we choose the same $ h_{\lambda} $, except the fact that here we use the computation with $ M_i^{-(1+ \beta )} $ in front the regular part of $ h_{\lambda} $.

Here also, we consider $ r_i \in (0, L_i) $ where $ L_i $ is the number of the blow-up analysis.

$$  L_i= \dfrac{1}{2i}R_i M_i^{(1+ \beta)/2}. $$

We argue by contradiction and we suppose that:

$$ M_i^{ \beta }D_i  \to +\infty . $$

Then, without loss of generality we can assume that:

$$ \dfrac{\nabla V_i(x_i)}{|\nabla V_i(x_i)|}\to e=(1,0,...0). $$

We use the Kelvin transform twice and around this point and around 0.

$$ h_{\lambda} (y)= \epsilon r_i^{-2} G_{\lambda}(y, \dfrac{e}{\delta})- \int_{\Sigma_{\lambda}} G_{\lambda}(y, \eta) Q_{\lambda}(\eta) d\eta. $$

with,

$$ Q_{\lambda}(\eta)=(V_{\delta}(\eta)-V_{\delta}(\eta^{\lambda})) (v_i^{\delta}(\eta^{\lambda}))^{3}. $$

And, by the same estimates, we have for $ \eta  \in A_1 $

$$ V_{\delta}(\eta)-V_{\delta}(\eta^{\lambda})  \geq M_i^{-1}D_i (\eta_1-\lambda) + o(1)M_i^{-1}|\eta^{\lambda}| , $$

and, we have for $ \eta \in A_2, |I_{\delta}(\eta)| \leq c_2M_i D_i^{1/\beta } $,

$$ | V_{\delta}(\eta)-V_{\delta}(\eta^{\lambda}) |  \leq CM_i^{-1}D_i (|I_{\delta}(\eta)| + |I_{\delta}(\eta^{\lambda})|), $$

and for $  M_i D_i^{1/ \beta } \leq |I_{\delta}(\eta)| \leq r_i $,

$$ | V_{\delta}(\eta)-V_{\delta}(\eta^{\lambda}) |  \leq M_i^{-1}D_i|I_{\delta}(\eta)| + M_i^{-(1+ \beta) }|I_{\delta}(\eta)|^{(1+ \beta )}, $$

By the same estimates, we have for $|I_{\delta}(\eta)| \leq r_i $ or  $ |y-e/\delta| \geq  c_3 r_i^{-1} $:

$$ h_{\lambda} (y)\simeq \epsilon r_i^{-2} G_{\lambda}(y, \dfrac{e}{\delta})+ c_4M_i^{- 1}D_i \dfrac{(y_1- \lambda)}{|y|^n} + o(1) M_i^{-1}D_i \dfrac{(y_1- \lambda)}{|y|^n} + o(1) M_i^{- (1+ \beta )}  G_{\lambda}(y, \dfrac{e}{\delta}). $$

with $ c_4 >0 $.

And, we have for some $ \lambda_0 \leq -2 $ and $ C_0 >0 $:

$$ v_i^{\delta}(y)- v_i^{\delta}(y^{\lambda_0}) \geq C_0\dfrac{y_1-\lambda_0}{(1+|y|)^n}, $$

for $ y_1 > \lambda_0 $.

\bigskip

By the same estimates as in Chen-Lin paper (we apply the lemma 2.1 of the second section),  and by our hypothesis on $ v_i $, we have:

$$ 0 < h_{\lambda} (y) < 2 \epsilon <  w_{\lambda}(y). $$

also, we have the same etimate on the boundary, $|I_{\delta}(\eta)| = r_i $ or  $ |y-e/\delta| =  c_5 r_i^{-1} $

\bigskip

\underbar {Step 2.2 conclusion : a linear equation perturbed by a term, and, the auxiliary function } 

\bigskip

Here also, we use the computations of Chen-Lin, and, we take the same auxiliary function $ h_{\lambda} $ (which correspond to this step), except the fact that here in front the regular part of this function we have  $ M_i^{-(1+ \beta )} $. 

\bigskip

Here also, we consider $ r_i \in (0, L_i) $ where $ L_i $ is the number of the blow-up analysis.

$$  L_i= \dfrac{1}{2i}R_i M_i^{(1+ \beta )/2}. $$

We set,

$$ v_i(z)=v_i^*(z+e), $$

where $ v_i^* $ is the blow-up function. Then, $ v_i $ has a local maximum near $ -e $.

$$ U_0(z)=U_0^*(z+e). $$

We have, for $  |y|\geq L_i'^{-1} $, $ L_i'=\dfrac{1}{2} R_iM_i, $

$$ \bar v_i(y)=\dfrac{1}{|y|^{n-2}} v_i \left ( \dfrac{y}{|y|^2} \right ). $$

$$ | V_i(\bar x_i+M_i^{-1}\dfrac{y}{ |y|^2} )-  V_i(\bar x_i) |  \leq M_i^{- (1+ \beta )}(1+|y|^{-1}).  $$

$$ x_i=\bar x_i+ M_i^{-1} e, $$

Then, for simplicity, we can assume that,  $ \bar v_i $ has a local maximum near $ e^*=(-1/2, 0, ...0) $.

Also, we have:

$$ | V_i(x_i+M_i^{-1}\dfrac{y}{ |y|^2} )-  V_i(x_i+M_i^{-1}\dfrac{y^{\lambda}}{ |y^{\lambda}|^2} ) |  \leq M_i^{-(1+ \beta )}(1+|y|^{-1}).  $$

$$ h_{\lambda} (y)\simeq \epsilon r_i^{-2} G_{\lambda}(y, 0)-\int_{\Sigma'_{\lambda}} G_{\lambda}(y, \eta) Q_{\lambda}(\eta) d\eta. $$

where, $ \Sigma'_{\lambda} = \Sigma_{\lambda} - \{ \eta, |\eta | \leq r_i^{-1} \}, $ and,

$$  Q_{\lambda}(\eta)= \left (V_i(x_i+M_i^{-1}\dfrac{y}{ |y|^2} )-  V_i(x_i+M_i^{-1}\dfrac{y^{\lambda}}{ |y^{\lambda}|^2} ) \right )(v_i(y^{\lambda}))^{3}. $$

we have by the same computations that:

$$ \int_{\Sigma'_{\lambda}} G_{\lambda}(y, \eta) Q_{\lambda}(\eta) d\eta \leq C M_i^{-(1+ \beta )}G_{\lambda}(y, 0) << \epsilon r_i^{-2} G_{\lambda}(y, 0). $$

By the same estimates as in Chen-Lin paper (we apply the lemma 2.1 of the second section),  and by our hypothesis on $ v_i $, we have:

$$ 0 < h_{\lambda} (y)  < 2 \epsilon <  w_{\lambda}(y). $$

also, we have the same estimate on the boundary, $ |y| = \dfrac{1}{r_i}  $. 

\bigskip 

\underbar {Proof of the theorem 2, $ n=5 $: }

\bigskip

To prove the theorem, we argue by contradiction and we assume that the $ (\sup )^{1/3} \times \inf $ tends to infinity.

\bigskip

\underbar {Step 1: blow-up analysis }

\bigskip

We want to prove that:

$$  {\tilde R}^3 (\sup_{B_{\tilde R}(0)} u)^{1/3} \times  \inf_{B_{3{\tilde R}}(0)} u \leq c= c(a, b, A, B), $$

If it is not the case, we have:

$$ {\tilde R_i}^3 ( \sup_{B_{\tilde R_i}(0)} u_i)^{1/3} \times  \inf_{B_{3 \tilde R_i}(0)} u_i = i^6 \to + \infty, $$

For positive solutions $ u_i >0 $ of the equation $ (E) $ and $ \tilde R_i \to 0. $

Thus,

$$ \dfrac{1}{i}{\tilde R_i}(\sup_{B_{\tilde R_i}(0)} u_i)^{2/3} \to + \infty, $$

and,

$$  \dfrac{1}{i}{\tilde R_i} [\sup_{B_{\tilde R_i}(0)} u_i]^{4/9} \to + \infty, $$

Let $ a_i $ such that:

$$ u_i(a_i) = \max_{B_{\tilde R_i}(0)} u_i, $$

We set, 

$$ s_i(x)=(\tilde R_i-|x-a_i|)^{9/4} u_i(x), $$

we have,

$$  s_i(\bar x_i)= \max_{B_{\tilde R_i}(a_i)} s_i \geq s_i(a_i)={\tilde R_i}^{9/4}\sup_{B_{R_i}(0)} u_i \to + \infty, $$

we set, 

$$ R_i= \dfrac{1}{2} (\tilde R_i-|\bar x_i-a_i|), $$

We have, for $ |x-\bar x_i|\leq \dfrac{R_i}{i} $, 

$$ \tilde R_i-|x-a_i|\geq \tilde R_i-|\bar x_i-a_i|-|x-a_i| \geq 2R_i-R_i=R_i$$

Thus,

$$ \dfrac{u_i(x)}{u_i(\bar x_i)} \leq \beta_i \leq 2^{9/4}. $$

with $ \beta_i \to 1. $

\bigskip

We set, 

 $$  M_i= u_i(\bar x_i),  \,\,\,v_i^*(y)=\dfrac{u_i(\bar x_i+M_i^{-2/3}y)}{u_i(\bar x_i)}, \,\,\, |y|\leq \dfrac{1}{i} R_iM_i^{4/9} = 2L_i. $$

And,

$$ \dfrac{1}{i^3}{R_i}^3 M_i^{1/3} \times  \inf_{B_{3 \tilde R_i}(0)} u_i \to + \infty, $$

By the elliptic estimates, $ v_i^* $ converge on each compact set of $ {\mathbb R}^5 $ to a function $ U_0^* >0 $ solution of :

$$ \begin {cases}

   -\Delta U_0^* = V(0) {U_0^*}^{7/3} \,\, \text{in} \,\, {\mathbb R}^5, \\
                
                U_0^*(0)=1= \max_{{\mathbb R}^5} U_0^*.

\end {cases} $$

For simplicity, we assume that $ 0 < V(0)= n(n-2) =15 $. By a result of Caffarelli-Gidas-Spruck, see [10], we have:

$$ U_0^*(y)=(1+|y|^2)^{-3/2}. $$

We set,

$$ v_i(y)=v_i^*(y+e), $$

where $ v_i^* $ is the blow-up function. Then, $ v_i $ has a local maximum near $ -e $.

$$ U_0(y)=U_0^*(y+e). $$

We want to prove that:

$$ \min_{ \{0 \leq |y| \leq r  \} } v_i^*  \leq (1+ \epsilon) U_0^*(r). $$

for $ 0 \leq r \leq L_i $, with $  L_i= \dfrac{1}{2i}R_i M_i^{4/9} $.

\bigskip

We assume that it is not true, then, there is a sequence of number $ r_i \in (0, L_i) $ and $ \epsilon >0 $, such that:

$$  \min_{ \{ 0 \leq |y| \leq r_i  \} } v_i^*  \geq (1+ \epsilon) U_0^*(r_i). $$

We have:

$$ r_i \to + \infty. $$

Thus , we have for $ r_i \in (0, L_i) $ :

$$  \min_{ \{ 0 \leq |y| \leq r_i  \} } v_i  \geq (1+ \epsilon) U_0(r_i). $$

Also, we can find a sequence of number $ l_i  \to + \infty $ such that:

$$ ||v_i^*-U_0||_{C^2(B_{l_i}(0))} \to 0. $$

Thus,

 $$ \min_{ \{ 0 \leq |y| \leq l_i  \} } v_i  \geq (1 - \epsilon/2) U_0(l_i). $$

\underbar {Step 2 : The Kelvin transform and the Moving-plane method }

\bigskip

\underbar {Step 2.1: a linear equation perturbed by a term, and, the auxiliary function } 

\bigskip

 {Step 2.1.1: $ D_i= |\nabla V_i(x_i)| \to 0.$} 

\bigskip

We have the same estimate as in the paper of Chen-Lin. We argue by contradiction. We consider $ r_i \in (0, L_i) $ where $ L_i $ is the number of the blow-up analysis.

$$  L_i= \dfrac{1}{2i}R_i M_i^{4/9}. $$

 We use the assumption that the sup times inf is not bounded to prove $ w_{\lambda} >h_{\lambda} $ in $ \Sigma_{\lambda} =\{ y, y_1 > \lambda \}, $ and on the boundary.
 
 \bigskip

The function $ v_i $ has a local maximum near $ -e $ and converge to $ U_0(y)=U_0^*(y+e) $ on each compact set of $ {\mathbb R}^5 $. $U_0 $ has a maximum at $ -e $.

We argue by contradiction and we suppose that:

$$ D_i= |\nabla V_i(x_i)| \not  \to 0. $$

Then, without loss of generality we can assume that:

$$ \nabla V_i(x_i) \to e=(1,0,...0). $$

Where $ x_i $ is :

$$ x_i=\bar x_i+ M_i^{-2/3} e, $$

with $ \bar x_i  $ is the local maximum in the blow-up analysis.

\bigskip

As in the paper of Chen-Lin, we use the Kelvin transform twice and we set (we take the same notations):

$$ I_{\delta}(y)=\dfrac{\frac{|y|}{|y|^2}-\delta e}{\left (|\frac{|y|}{|y|^2}-\delta e| \right )^2}, $$

$$ v_i^{\delta}(y)=\dfrac {v_i(I_{\delta}(y))}{|y|^{n-2}|y-e/\delta|^{n-2}}, $$

and, 

$$  V_{\delta}(y) =V_i(x_i+ M_i^{-2/3}I_{\delta}(y)). $$

$$ U_{\delta}(y)=\dfrac {U_0(I_{\delta}(y))}{|y|^{n-2}|y-e/\delta|^{n-2}}. $$

Then, $ U_{\delta} $ has a local maximum near $ e_{\delta} \to -e $ when $ \delta \to 0 $. The function $ v_i^{\delta} $ has a local maximum near $ -e $.

We want to prove by the application of the maximum principle and the Hopf lemma that near $ e_{\delta} $ we have not a local maximum, which is a contradiction.

\bigskip

We set on $ \Sigma'_{\lambda} = \Sigma_{\lambda}  - \{ y, |y- \frac{e}{\delta}| \leq \frac{c_0}{r_i} \} \simeq \Sigma_{\lambda}  - \{ y, |I_{\delta}(y)| \geq  r_i \} $:

$$ h_{\lambda} (y)= - \int_{\Sigma_{\lambda}} G_{\lambda}(y, \eta) Q_{\lambda}(\eta) d\eta. $$
 
 with,

$$ Q_{\lambda}(\eta)=(V_{\delta}(\eta)-V_{\delta}(\eta^{\lambda})) (v_i^{\delta}(\eta^{\lambda}))^{(n+2)/(n-2)}. $$

And, by the same estimates, we have for $ \eta  \in A_1=\{ \eta, |\eta | \leq R=\epsilon_0/\delta \},$

$$ V_{\delta}(\eta)-V_{\delta}(\eta^{\lambda})  \geq M_i^{-2/3}(\eta_1-\lambda) + o(1)M_i^{-2/3}|\eta^{\lambda}| , $$

and, we have for $ \eta \in A_2 =\Sigma_{\lambda}-A_1 $:

$$ | V_{\delta}(\eta)-V_{\delta}(\eta^{\lambda}) |  \leq CM_i^{-2/3}(|I_{\delta}(\eta)| + |I_{\delta}(\eta^{\lambda})|), $$

And, we have for some $ \lambda_0 \leq -2 $ and $ C_0 >0 $:

$$ v_i^{\delta}(y)- v_i^{\delta}(y^{\lambda_0}) \geq C_0\dfrac{y_1-\lambda_0}{(1+|y|)^n}, $$

for $ y_1 > \lambda_0 $.

\bigskip

By the same estimates,  and by our hypothesis on $ v_i $, we have, for $ c_1 >0 $:

$$ 0 < h_{\lambda} (y) < 2 \epsilon <  w_{\lambda}(y). $$

also, we have the same etimate on the boundary, $|I_{\delta}(\eta)| = r_i $ or  $ |y-e/\delta| =  c_2 r_i^{-1} $

\bigskip

Step 2.1.1: $ |\nabla V_i(x_i)|[u_i(x_i)]^{2/3} \leq C $

\bigskip

Here, also, we argue by contradiction. We use the same computation as in Chen-Lin paper, we take $ \alpha = 2 $ and we choose the same $ h_{\lambda} $, except the fact that here we use the computation with $ M_i^{-4/3} $ in front the regular part of $ h_{\lambda} $.

Here also, we consider $ r_i \in (0, L_i) $ where $ L_i $ is the number of the blow-up analysis.

$$  L_i= \dfrac{1}{2i}R_i M_i^{4/9}. $$

We argue by contradiction and we suppose that:

$$ M_i^{2/3}D_i  \to +\infty . $$

Then, without loss of generality we can assume that:

$$ \dfrac{\nabla V_i(x_i)}{|\nabla V_i(x_i)|}\to e=(1,0,...0). $$

We use the Kelvin transform twice and around this point and around 0.

$$ h_{\lambda} (y)= \epsilon r_i^{-3} G_{\lambda}(y, \dfrac{e}{\delta})- \int_{\Sigma_{\lambda}} G_{\lambda}(y, \eta) Q_{\lambda}(\eta) d\eta. $$

with,

$$ Q_{\lambda}(\eta)=(V_{\delta}(\eta)-V_{\delta}(\eta^{\lambda})) (v_i^{\delta}(\eta^{\lambda}))^{(n+2)/(n-2)}. $$

And, by the same estimates, we have for $ \eta  \in A_1 $

$$ V_{\delta}(\eta)-V_{\delta}(\eta^{\lambda})  \geq M_i^{-2/3}D_i (\eta-\lambda) + o(1)M_i^{-2/3}|\eta^{\lambda}| , $$

and, we have for $ \eta \in A_2, |I_{\delta}(\eta)| \leq c_2M_i^{2/3}D_i $,

$$ | V_{\delta}(\eta)-V_{\delta}(\eta^{\lambda}) |  \leq CM_i^{-2/3}D_i (|I_{\delta}(\eta)| + |I_{\delta}(\eta^{\lambda})|), $$

and for $  M_i^{2/3}D_i \leq |I_{\delta}(\eta)| \leq r_i $,

$$ | V_{\delta}(\eta)-V_{\delta}(\eta^{\lambda}) |  \leq M_i^{-2/3}D_i|I_{\delta}(\eta)| + M_i^{-4/3}|I_{\delta}(\eta)|^2, $$

By the same estimates, we have for $|I_{\delta}(\eta)| \leq r_i $ or  $ |y-e/\delta| \geq  c_3 r_i^{-1} $:

$$ h_{\lambda} (y)\simeq \epsilon r_i^{-3} G_{\lambda}(y, \dfrac{e}{\delta})+ c_4M_i^{-2/3}D_i \dfrac{(y_1- \lambda)}{|y|^n} + o(1) M_i^{-2/3}D_i \dfrac{(y_1- \lambda)}{|y|^n} + o(1) M_i^{-4/3}  G_{\lambda}(y, \dfrac{e}{\delta}). $$

with $ c_4 >0 $.

And, we have for some $ \lambda_0 \leq -2 $ and $ C_0 >0 $:

$$ v_i^{\delta}(y)- v_i^{\delta}(y^{\lambda_0}) \geq C_0\dfrac{y_1-\lambda_0}{(1+|y|)^n}, $$

for $ y_1 > \lambda_0 $.

\bigskip

By the same estimates as in Chen-Lin paper (we apply the lemma 2.1 of the second section),  and by our hypothesis on $ v_i $, we have:

$$ 0 < h_{\lambda} (y) < 2 \epsilon <  w_{\lambda}(y). $$

also, we have the same etimate on the boundary, $|I_{\delta}(\eta)| = r_i $ or  $ |y-e/\delta| =  c_5 r_i^{-1} $:

\bigskip

\underbar {Step 2.2 conclusion : a linear equation perturbed by a term, and, the auxiliary function } 

\bigskip

Here also, we use the computations of Chen-Lin, and, we take the same auxiliary function $ h_{\lambda} $ (which correspond to this step), except the fact that here in front the regular part of this function we have  $ M_i^{-4/3} $. 

\bigskip

Here also, we consider $ r_i \in (0, L_i) $ where $ L_i $ is the number of the blow-up analysis.

$$  L_i= \dfrac{1}{2i}R_i M_i^{4/9}. $$

We set,

$$ v_i(z)=v_i^*(z+e), $$

where $ v_i^* $ is the blow-up function. Then, $ v_i $ has a local maximum near $ -e $.

$$ U_0(z)=U_0^*(z+e). $$

We have, for $  |y|\geq L_i^{-1} $, $ L_i=\dfrac{1}{2} R_iM_i^{2/3}, $

$$ \bar v_i(y)=\dfrac{1}{|y|^{n-2}} v_i \left ( \dfrac{y}{|y|^2} \right ). $$

$$ | V_i(\bar x_i+M_i^{-2/3}\dfrac{y}{ |y|^2} )-  V_i(\bar x_i) |  \leq M_i^{-4/3}(1+|y|^{-2}).  $$

$$ x_i=\bar x_i+ M_i^{-2/3} e, $$

Then, for simplicity, we can assume that,  $ \bar v_i $ has a local maximum near $ e^*=(-1/2, 0, ...0) $.

Also, we have:

$$ | V_i(x_i+M_i^{-2/3}\dfrac{y}{ |y|^2} )-  V_i(x_i+M_i^{-2/3}\dfrac{y^{\lambda}}{ |y^{\lambda}|^2} ) |  \leq M_i^{-4/3}(1+|y|^{-2}).  $$

$$ h_{\lambda} (y)\simeq \epsilon r_i^{-3} G_{\lambda}(y, 0)-\int_{\Sigma'_{\lambda}} G_{\lambda}(y, \eta) Q_{\lambda}(\eta) d\eta. $$

where, $ \Sigma'_{\lambda} = \Sigma_{\lambda} - \{ \eta, |\eta | \leq r_i^{-1} \}, $ and,

$$  Q_{\lambda}(\eta)= \left (V_i(x_i+M_i^{-2/3}\dfrac{y}{ |y|^2} )-  V_i(x_i+M_i^{-2/3}\dfrac{y^{\lambda}}{ |y^{\lambda}|^2} ) \right )(v_i(y^{\lambda}))^{\frac{n+2}{n-2}}. $$

we have by the same computations that:

$$ \int_{\Sigma'_{\lambda}} G_{\lambda}(y, \eta) Q_{\lambda}(\eta) d\eta \leq C M_i^{-4/3}G_{\lambda}(y, 0) << \epsilon r_i^{-3} G_{\lambda}(y, 0). $$

By the same estimates as in Chen-Lin paper (we apply the lemma 2.1 of the second section),  and by our hypothesis on $ v_i $, we have:

$$ 0 < h_{\lambda} (y)  < 2 \epsilon <  w_{\lambda}(y). $$

also, we have the same estimate on the boundary, $ |\eta| = \dfrac{1}{r_i}  $.

\end{document}